\documentclass[AMA,STIX1COL]{WileyNJD-v2}

\articletype{Article Type}%

\received{--}
\revised{--}
\accepted{--}

\usepackage{color}
\usepackage{soul}

\begin{document}

\title{An implicit locking-free B-spline Material Point Method for large strain geotechnical modelling}

\author[1]{Mian Xie}

\author[2]{Pedro Navas}

\author[1]{Susana L\'opez-Querol*}

\authormark{Xie, M. \textsc{et al}}

\address[1]{\orgdiv{University College London.}, \orgname{Department of Civil, Environmental and Geomatic Engineering.}, \orgaddress{\state{London.}, \country{UK}}}

\address[2]{\orgdiv{Universidad Polit\'ecnica de Madrid.}, \orgname{ETSI Caminos, Canales y Puertos.}, \orgaddress{\state{Madrid.}, \country{Spain}}}

\corres{*Susana L\'opez-Querol. \email{s.lopez-querol@ucl.ac.uk}}


\abstract[Summary]{The Material Point Method (MPM) has drawn great attention in the numerical modelling of large deformation, geotechnical problems. The popularity of MPM is mainly because its formulation shares significant similarities with the finite element method. In MPM, the iteration points can move independently from the mesh, allowing for the resolution of large deformation problems. However, because of this, the original MPM formulation suffers from the well-known cell-crossing noise and volumetric-locking instabilities, resulting in a strongly oscillated stress field. A novel implicit locking-free B-spline MPM that controls stress oscillations to a negligible level is proposed in this paper. A novel, but very straightforward B-spline shape function implementation procedure, avoids the need for a complex material point searching algorithm, providing seamless transformation from the original MPM to this robust B-spline MPM, aiming at modelling large-strain geotechnical problems. The newly proposed volumetric locking mitigation strategy is also very easy to implement, which facilitates the reproducibility of this research. The proposed method is validated against three numerical studies: granular column collapse experiment, slope failure and footing with large penetration. The proposed numerical method agrees well with experiments reported in the literature and previous numerical studies. Also, these numerical examples show that the proposed method provides a more prominent stress field than other available methodologies.}


\jnlcitation{\cname{%
\author{M. Xie},
\author{P. Navas} and 
\author{S. L\'opez-Querol}} (\cyear{2023}), 
\ctitle{An implicit locking-free B-spline Material Point Method for large strain geotechnical modelling}, \cjournal{I.J. Num. Anal. Meth. Geomechanics}, \cvol{2023;00:1--14}.}

\maketitle


\section{Introduction and literature review}
\label{sec:1}
From the design point of view, civil engineers are often interested in small deformation problems because it is feasible to undertake the analysis of structures that are not likely to experience high strains. The Finite Element Method (FEM) and the Finite Difference Method (FDM) are the most commonly used numerical approaches in the civil engineering industry due to their accuracy for small deformation problems and the ability to describe the history-dependent material behaviour. However, geotechnical engineering involves numerous large deformation problems, such as landslides and penetration of piles during their installation, among others. 

Soga et al.~\cite{Soga2016} comprehensively reviewed the available numerical methods for large deformation in geotechnical modelling. The Material Point Method (MPM) is recommended due to the following reasons: (a) it is suitable for large-scale and large deformation problems; (b) its computational efficiency; (c) it is based on a similar mathematical formulation as for the FEM; (d) its ability to deal with the history-dependent soil constitutive models; and (e) the simplicity for applying essential boundary conditions~\cite{Soga2016}. { Some other comprehensive reviews of numerical modelling for large deformation in geotechnics can be found in ~\cite{Wang2015, Augarde2021, Qin2022, Shan2022}. }

In MPM, the material domain is discretised in material points that can be treated as integration points in FEM~\cite{Coombs2018}. To solve the distorted mesh issue in FEM, the mesh is set to its original position in MPM at each time step, and the material points carry all the permanent information (e.g. stress, displacement and velocity). As a result, the material points (integration points) can move independently from the mesh~\cite{Coombs2018}. In fact, the MPM will yield to FEM if the mesh is not set to its original position after every time step~\cite{Guilkey2003}. 

Aimed for computational efficiency, Sulsky et al.~\cite{Sulsky1994} use linear FEM shape functions (hat functions) in the original MPM. The use of these shape functions in MPM causes well-known cell crossing instability when a material point crosses the cell boundary. The Generalised Interpolation MPM (GIMPM) developed by Bardenhagen and Kober~\cite{Bardenhagen2004} is one of the most commonly used MPM methods in the literature due to its efficiency. However, it cannot fully eliminate the cell crossing instability~\cite{Zhang2016} since for that it is necessary to utilise a higher-order shape function that covers multiple cells.

The available higher-order MPMs are the B-spline MPM (BSMPM)~\cite{Steffen2008} and the Local Maximum-Entropy MPM~\cite{Molinos2021a}. Both methods can entirely eliminate the cell crossing noise. BSMPM is chosen for this research due to its better efficiency. The latter approach requires a Newton-Raphson loop to construct the Local Maximum-Entropy shape function~\cite{Arroyo2006, Molinos2021b, Molinos2021c}. However, the construction of the B-spline shape function still requires much more computational effort than the original MPM and the GIMPM shape function. Additionally, the conventional implementation of a higher-order shape function requires a particle-searching algorithm to know where the material points and their neighbours are located~\cite{Molinos2021c}. Therefore, it is more appropriate to use an implicit iteration algorithm rather than an explicit one for the MPM with higher-order shape functions, since an implicit solver allows relatively larger time steps and the shape function only needs to be constructed once for each step. Also, due to this feature, the implicit solver is more preferable for relatively slow simulations, especially for the quasi-static analysis. {Apart from utilising the higher-order shape function, the Moving Least Squares (MLS) shape function is also usually used to mitigate the cell crossing noise~\cite{Hu2018, Tran2020, Song2021}. In the current literature, the MLS shape function has only been used in the explicit MPM solver, and the implicit MLS-MPM still needs to be studied.}

Volumetric locking is the other source of instability in MPM besides the cell crossing noise. The volumetric locking will result in over-stiff behaviour and non-physical stress oscillations. In FEM, the volumetric locking instability {can be easily overcome} by a reduced quadrature scheme~\cite{Bower2009}. However, utilising the reduced quadrature scheme is not feasible in MPM for two reasons: (a) a sufficient number of material points per cell is required to ensure the desired quadrature accuracy; (b) it is difficult to limit the number of material points per cell since they can move freely between them~\cite{Coombs2018}. The available volumetric locking mitigation strategies~\cite{Coombs2018, Love2006, Mast2012, Kularathna2017, Zhang2017, Iaconeta2019, Bisht2021, Wang2021, Zheng2022} are mainly suitable for the MPM with linear shape functions (e.g. original MPM and GIMPM). The BSMPM shows less volumetric locking than the MPMs with linear shape functions. However, fully overcoming the volumetric locking for BSMPM is challenging since its shape functions are associated with adjacent cells~\cite{Zhao2022}. Navas et al.~\cite{Navas2018} mitigated the volumetric locking in the Optimal Transportation Meshfree (OTM) method, a numerical method based on Local Maximum-Entropy shape functions. This approach relies on the triangular mesh, however, the rectangular mesh is a more popular choice in the MPM literature. To the best of our knowledge, there are currently only three available pieces of research~\cite{Telikicherla2022, Zhao2022, Sugai2023} related to the mitigation of volumetric locking for a higher-order MPM. Two of them~\cite{Telikicherla2022, Sugai2023} are based on the F-bar projection method proposed by Elguedj et al.~\cite{Elguedj2008}, which was originally developed for B-spline FEM. However, the application of this method is not straightforward due to the introduction of an additional background grid with low-order shape functions~\cite{Zhao2022}. {Also, it is hard to apply to other types of MPM (i.e. GIMPM) because of this feature.} The method proposed by Zhao et al.~\cite{Zhao2022} is simple enough, but its performance has not been well tested by a wide range of large deformation geotechnical problems. 

This paper proposes a total locking-free, implicit BSMPM implementation to overcome the previously mentioned research gaps. This newly proposed method is simple, as does not require the implementation of any sophisticated particle search algorithm. The application of the newly proposed volumetric locking migration strategy is also very straightforward and can be easily implemented. Additionally, the proposed method has been validated with various large deformation geotechnical problems.

This paper is structured as follows: Section \ref{sec:2} explains the proposed numerical model in detail, including the construction of B-spline shape functions, large strain formulation of the constitutive model, volumetric locking migration strategy and the proposed implicit BSMPM algorithm. Then, we will present some verification numerical examples for large deformation geotechnical problems in Section \ref{sec:3}: the granular column collapse, slope failures and the large penetration of strip footing problem. Finally, we present the conclusions of this research in Section \ref{sec:4}.

\section{Description of the numerical model}
\label{sec:2}
\subsection{BSMPM shape function: a new efficient searching algorithm}
\label{subsec:21}
In MPM, the kinematic quantities (e.g. displacement, velocity and acceleration) at any material point can be approximated by interpolating the nodal quantities of the related grid using the shape functions. For example, the incremental displacement of a material point ($\Delta \pmb{u}_{p}$) can be approximated from the incremental nodal displacement of related grid nodes ($\pmb{u}_{I}$) by
\begin{equation}\label{eq1}
	\Delta\pmb{u}_{p} = \sum_{I=1}^{{n}_{n}} {N}_{Ip}\Delta\pmb{u}_{I},
\end{equation}
where ${n}_{n}$ is the number of grid nodes that influences the material point. In this research, the subscript $p$ means the variables are related to the material point. Similarly, the subscript $I$ denotes the variables are related to the grid nodes. ${N}_{Ip}$ is the shape function associated with node $I$ evaluated at the position of material point $p$. In the original MPM, the employed shape function is the traditional hat function; however, different shape functions are used for improved MPMs, such as the B-spline function. Also, the total nodal displacements are equivalent to the incremental nodal displacements in MPM once a brand-new grid has been used for each time step.

The communication between nodes and material points is not one-way only. The information can also be mapped from the material points to the related node. For example, the {diagonal components of lumped} nodal mass ($m_I$) can be calculated by taking into account the contribution of all the material points related to this node:
\begin{equation}\label{eq2}
	m_I = \sum_{p=1}^{{n}_{p}} {N}_{Ip} m_p,
\end{equation}
where ${n}_{p}$ is the number of material points associated with node $I$ and $m_{p}$ is the mass of the material point $p$.

BSMPM follows the exactly same computational procedures as the original MPM except for utilising a different shape function ${N}_{Ip}$~\cite{Motlagh2017}. Like the other MPM shape functions, we can first construct the B-spline shape functions in 1-D and then convolute them in multi-dimensions~\cite{Gan2018, Tran2019, Wobbes2019, Sun2022}. A 1-D BSMPM shape function is usually constructed based on a knot vector, $\Xi$, which contains a series of knots with non-decreasing order for their spatial coordinates:
\begin{equation}\label{eq3}
	\Xi = \{ {x}_{1}, {x}_{2}, \cdots, {x}_{n+q}, {x}_{n+q+1} \},
\end{equation}
where $n$ is the total number of B-spline shape functions; $q$ denotes the polynomial order and $ {x}_{i} (i=1,2,\cdots,n+q+1) $ is the $i^{th}$ knot~\cite{Wobbes2019}. These knots can be treated as grid nodes in the FEM or MPM. However, in the case of a boundary grid node, the corresponding boundary knot needs to be repeated $q$ times. As we can see from the definition of the knot vector, at least $q+2$ knots are required to construct a B-spline shape function. In this research, we use a quadratic ($q=2$) B-spline shape function; therefore, at least four knots are needed to construct a B-spline shape function. In the case of a boundary grid, only two grid nodes are available. To have four knots, its boundary knot needs to be repeated twice.

\begin{figure}
	\begin{center}
		\includegraphics[width=1.0\textwidth]{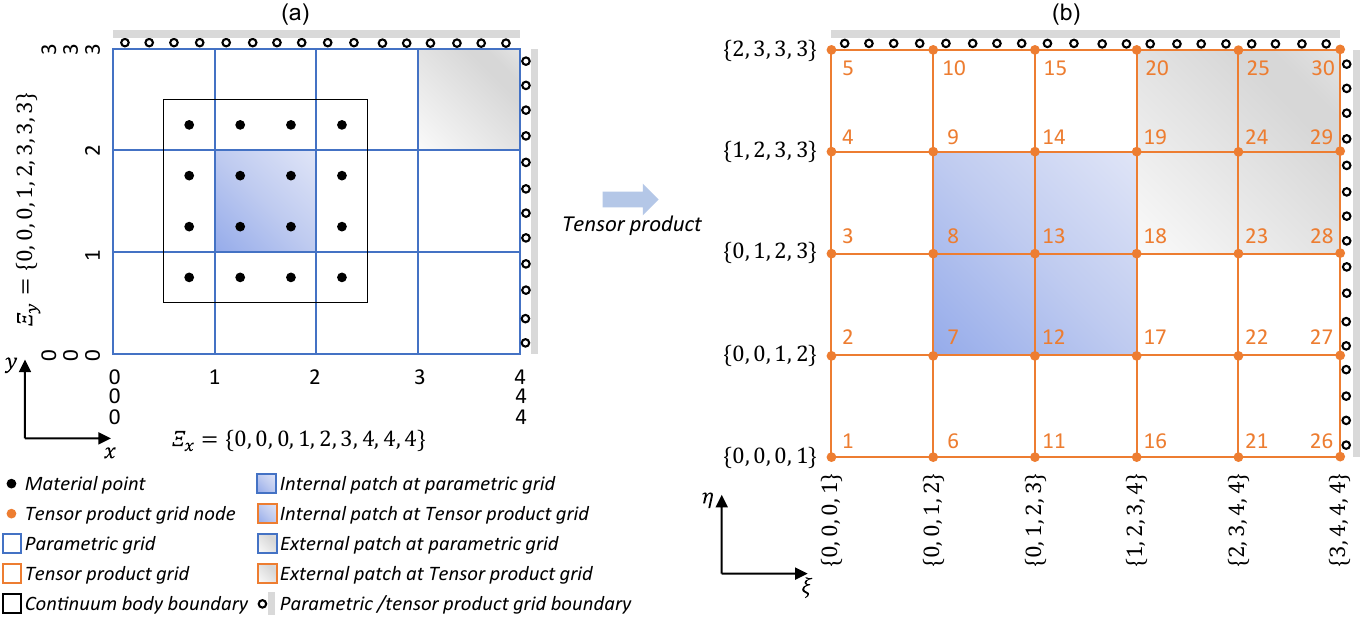}
		\caption{Graphical illustration of 2-D quadratic B-spline MPM discretisation: (a) parametric grid, (b) tensor product grid}
		\label{fig1}
	\end{center}
\end{figure}

Fig.\ref{fig1} graphically illustrates the 2-D discretisation of quadratic BSMPM. Fig.\ref{fig1}(a) shows a discretisation of a rectangular continuum body using sixteen material points. This continuum body moves in a domain (parametric grid) uniformly discretised into twelve patches (four patches in the x-direction and three patches in the y-direction) with a constant unit spacing. The entire parametric grid can be represented by two open global knot vectors: ${\Xi}_{x}=\{0,0,0,1,2,3,4,4,4\}$ and ${\Xi}_{y}=\{0,0,0,1,2,3,3,3\}$ in x- and y-directions, respectively. The knot vectors are open because their first and last knot appear $q+1$ times. As shown in Fig.\ref{fig1}(a), these knots appear three times for a quadratic B-spline, and these repeated knots are the boundary knots. Fig.\ref{fig1}(b) shows a 2-D tensor product grid formed by the tensor product of the knot vectors in the x- and y-directions. In BSMPM, the information is mapped from the material points to the tensor product grid nodes~\cite{Sun2022}. Then, the partial differential equations are solved in the tensor product grid. Finally, all the permanent information is mapped back to the material points, and the tensor product grid is set to its original position. As shown in Fig.\ref{fig1}, the parametric grid and the tensor product grid only coincide at the boundary grid nodes where the Dirichlet boundary conditions apply~\cite{Wobbes2019}. 

We use an internal patch (blue hatched in Fig.\ref{fig1}) which contains four material points, as an example to illustrate the construction of the BSMPM shape function. To construct shape functions for each material point inside this patch, we need to obtain sub-knot vectors from the global knot vectors ${\Xi}_{x}$ and ${\Xi}_{y}$. As shown in Fig.\ref{fig1}(a), this internal patch sits on the second span of the parametric grid in both directions. Therefore, we need to obtain sub-knot vectors from the second knot in the global knot vectors ${\Xi}_{x}$ and ${\Xi}_{y}$. For each patch, there are $q+1$ shape functions in each direction, and $q+2$ knots are required to construct a shape function. As a result, the sub knot vectors $\{0, 0, 1, 2\} \{0, 1, 2, 3\} \{1, 2, 3, 4\}$ in the x-direction and $\{0, 0, 1, 2\} \{0, 1, 2, 3\} \{1, 2, 3, 3\}$ in the y-direction are obtained by reading four knots three times one by one from the second knot in the global knot vectors.

Given the coordinate of a material point ($x$) and a knot vector ($\Xi = \{ {x}_{1}, {x}_{2}, \cdots, {x}_{n+q}, {x}_{n+q+1} \}$), we can use the Cox-de Boor formula~\cite{deBoor2001} to determine the shape functions recursively starting from polynomial order zero ($q=0$):
\begin{equation}\label{eq4}
	{N}_{i,0}(x) = \Big\{
	\begin{matrix}
		1 \\ 0
	\end{matrix}
	\,\,\,\,
	\begin{matrix}
		if\, {x}_{i}\leq{x}<{x}_{i+1} \\ otherwise
	\end{matrix}.
\end{equation}
In the case of polynomial order larger than zeros ($q>0$):
\begin{equation}\label{eq5}
	{N}_{i,q}(x) = \frac{{x}-{x}_{i}}{{x}_{i+q}-x_{i}} {N}_{i,q-1}(x) + \frac{{x}_{i+q+1}-{x}}{{x}_{i+q+1}-x_{i+1}} {N}_{i+1,q-1}(x).
\end{equation}
The shape function derivatives can be obtained using the following:
\begin{equation}\label{eq6}
	\frac{\text{d}{N}_{i,q}(x)}{\text{d}x} = \frac{q}{{x}_{i+q}-x_{i}} {N}_{i,q-1}(x) - \frac{q}{{x}_{i+q+1}-x_{i+1}} {N}_{i+1,q-1}(x).
\end{equation}
Note that the above Cox-de Boor recursion formula assumes $0/0=0$.

As a result, there are three shape functions (i.e. $i=3$) for each material point in each direction. After the tensor product, there are nine (three time three) shape functions evaluated at the tensor product grid nodes for each material point. In other words, one material point is governed by nine tensor product grid nodes. These four material points in the hatched internal patch are governed by the same tensor product grid nodes (i.e. grid nodes number $7, 8, 9, 12, 13, 14, 17, 18$ and $19$), as shown in Fig.\ref{fig1}(b). Therefore, their material point ($p$) to tensor product grid nodes ($N$) connectivity matrix ($p2N$ matrix) is $[7, 8, 9, 12, 13, 14, 17, 18, 19]$. Similarly, if a material point moves to the hatched external patch, its $p2N$ matrix will be $[18, 19, 20, 23, 24, 25, 28, 29, 30]$, as shown in Fig.\ref{fig1}(b). This $p2N$ matrix is crucial in MPM during the communication between the material points and their related grid nodes (in BSMPM, they are tensor product grid nodes). Summations over the grid nodes and material points require this $p2N$ matrix to know the relationship between the material points and grid nodes, as shown in Eq.~(\ref{eq1}) and ~(\ref{eq2}). The conventional method constructs the BSMPM shape functions patch by patch. Utilising a particle search algorithm, we can determine the material point ($p$) to patch/element ($e$) connectivity matrix ($p2e$). Since the patch ($e$) to tensor product grid nodes ($N$) connectivity matrix ($e2N$) is an inherent property of a tensor product grid, we can obtain the $p2N$ matrix through $p2e$ and $e2N$ (i.e. material point to patch to tensor product grid nodes).

This research proposes a new implementation to avoid the searching algorithm. The basic idea is constructing shape functions for the entire parametric grid instead of patches. As shown in Fig.\ref{fig1}(a), for each material point, we can construct its B-spline shape functions for the entire background mesh instead of each patch using only two global knot vectors ${\Xi}_{x}$ and ${\Xi}_{y}$. Using the Cox-de Boor formula, we can obtain the shape functions evaluated at every tensor product grid node for each material point. However, these shape functions are non-zero only at the patch where the interested material point locates. Because of this property, we can eliminate the zero shape functions to obtain the exact same shape functions as constructing them patch by patch. Also, we can find and record the locations (indexes) of these non-zero shape functions. These indexes are the components of the $p2N$ matrix for this material point. Finally, this method yields the same shape functions and $p2N$ matrix as the conventional patch-by-patch method with a particle searching algorithm. {The algorithm of the proposed method for a 2-D problem is shown in Appendix \ref{ap:1}.}

{A quadratic B-spline is used in this research because of its efficiency and capability of fully overcoming the cell crossing noise~\cite{Steffen2008}. Although a cubic B-spline shape function results in a more accurate result~\cite{Sun2022}, it is more cumbersome compared with the quadratic one because of the larger matrix size. The construction of a cubic B-spline shape function follows the same recursion procedure (Eq.~(\ref{eq4}) to (\ref{eq6})) and the proposed methodology can be easily extended to the employment of such a higher order shape function.}

{In the computer graphics community~\cite{Stomakhin2013, Klar2016}, the explicitly defined B-spline shape function is usually used. This type of B-spline shape function is derived recursively assuming a uniform grid space and internal zero-centred material point~\cite{Steffen2008}. The employment of this type of B-spline shape function is very efficient because the recursion procedure is avoided during construction. However, implementing this explicitly defined B-spline shape function requires ghost nodes (i.e. additional nodes) outside the boundary to ensure the partition of unity~\cite{Steffen2008b}. The presence of these ghost nodes introduces massive numerical errors around the boundaries. In the case of MPM application in computer graphics~\cite{Stomakhin2013, Klar2016}, the material body is usually not placed against the boundaries. However, the massive numerical errors around the boundaries are not acceptable for an engineering problem because the symmetric boundary is frequently used to reduce the computational effort, for example, the strip footing problem in geotechnical engineering (Section~\ref{subsec:33}). Therefore, it is more appropriate to use a knot vector to construct the B-spline shape functions recursively (Eq.~(\ref{eq4}) to (\ref{eq6})) to avoid introducing ghost nodes.}

\subsection{Large strain constitutive framework}
\label{subsec:22}
In MPM, the constitutive model is evaluated at the material point level. Therefore, for the notations in this section, we omit the subscript $p$ for clarity purposes. 

In the large strain constitutive model, the left elastic Cauchy-Green strain ${\pmb{b}}^e$ is stored as a state variable, and it is updated every time step from $t$ to $t+\text{d}t$~\cite{deSouza2008}:
\begin{equation}\label{eq7}
	\pmb{b}^{e}(t+\text{d}t) = \Delta\pmb{F}\, \pmb{b}^e (t)\, \Delta\pmb{F}^T,
\end{equation}
{ where $\Delta\pmb{F}$ is the deformation gradient increment. The deformation gradient, $\pmb{F}$, can be written as:}

\begin{equation}\label{eq7b}
	\pmb{F} = \frac{\partial \pmb{x}}{\partial \pmb{X}} = \frac{\partial \left( \pmb{X}+\pmb{u} \right)}{\partial \pmb{X}} = \pmb{I} + \frac{\partial \pmb{u}}{\partial \pmb{X}}.
\end{equation}

One of the most commonly used strain measurements for the large strain formulation is logarithmic (also called Hencky strain tensor), and, therefore, the elastic logarithmic strain $\pmb{\epsilon}^e$ is defined by~\cite{deSouza2008}:
\begin{equation}\label{eq8}
	\pmb{\epsilon}^{e} = \frac{1}{2}\text{ln}\: \pmb{b}^{e}.
\end{equation}
Note that the evaluation of Eq.~(\ref{eq8}) requires polar decomposition since ${\pmb{b}}^e$ is a matrix. The corresponding stress measure is the Kirchhoff stress $\pmb{\tau}$ defined by
\begin{equation}\label{eq9}
	\pmb{\tau} = \pmb{D}^{e} \pmb{\epsilon}^{e},
\end{equation}
where $\pmb{D}^{e}$ is the linear elastic consistent tangent in Voigt notation. In the MPM formulation, the Cauchy stress $\pmb{\sigma}$ is still used~\cite{Coombs2020}. The Kirchhoff stress $\pmb{\tau}$ has only been used in the constitutive model. The Cauchy stress $\pmb{\sigma}$ can be recovered from the Kirchhoff stress $\pmb{\tau}$ by
\begin{equation}\label{eq10}
	\pmb{\sigma} = \frac{\pmb{\tau}}{\text{det}(\pmb{F})} = \frac{\pmb{\tau}}{J},
\end{equation}
where $J = \text{det}(\pmb{F})$ is the Jacobian of the deformation gradient.

With a continuum-based method, such as the Finite Element Method (FEM) or the Material Point Method (MPM), the stress needs to be updated along with the incremental displacement. Usually, the elastic predictor return mapping algorithm is used to update the stress in an elastoplastic constitutive model. Utilising the logarithmic strain and Kirchhoff stress allows us to use the same elastic predictor return mapping algorithm derived from the small strain theory~\cite{deSouza2008}. In this research, we use the Mohr-Coulomb constitutive model as {documented} by de Souza Neto et al.~\cite{deSouza2008}, who also report the return mapping algorithm in detail. 

\subsection{Overcoming volumetric locking: a new F-bar approach}
\label{subsec:23}
This research employs the F-bar method to overcome volumetric locking. This method is originally proposed by de Souza Neto et al.~\cite{deSouza1996} to avoid the volumetric locking in FEM. Coombs et al.~\cite{Coombs2018} summarised the advantages of the F-bar method as: (a) it does not introduce any additional unknown or unphysical parameter; (b) it is possible to apply to any type of constitutive model and MPM; (c) it does not require additional material points to capture the volumetric behaviour. Same as the previous section, the subscript $p$ is omitted for the notations in this section. The notation with no subsection means this value is evaluated at the material point level.

The idea of the F-bar method is based on the multiplicative split of the deformation gradient $\pmb{F}$
\begin{equation}\label{eq12}
	\pmb{F} = \pmb{F}_d\pmb{F}_v,
\end{equation}
where
\begin{equation}\label{eq13}
	\pmb{F}_d = \text{det}(\pmb{F})^{-1/dof}\pmb{F}  \;\;\;\; and \;\;\;\; \pmb{F}_v = \text{det}(\pmb{F})^{1/dof}\pmb{I}
\end{equation}
are the deviatoric and volumetric parts of the deformation gradient, respectively~\cite{deSouza1996}. $\pmb{I}$ is the identity matrix. $dof$ is the spatial dimensions of the problem (e.g. $dof=2$ for a 2-D problem).

In the F-bar method, the volumetric locking is resolved by replacing the deformation gradient $\pmb{F}$ in the constitutive model with a modified deformation gradient $\pmb{\Bar{F}}$~\cite{deSouza1996}. In this modified deformation gradient, the volumetric part $\pmb{F}_v$ of $\pmb{\Bar{F}}$ is replaced by $(\pmb{F}_0)_v$, a volumetric part of $\pmb{\Bar{F}}$ that has less volumetric constraint~\cite{deSouza1996}
\begin{equation}\label{eq14}
	\pmb{\Bar{F}} = \pmb{F}_d (\pmb{F}_0)_v = [\frac{\text{det}(\pmb{F}_0)}{\text{det}(\pmb{F})}]^{1/dof} \pmb{F}.
\end{equation}
In this research, we use the incremental version of the F-bar method since our constitutive model is driven by the incremental deformation gradient shown in Eq.~(\ref{eq7}). Therefore, the incremental version of the F-bar method can be obtained by rewriting Eq.~(\ref{eq14})~\cite{Coombs2018}:
\begin{equation}\label{eq15}
	\Delta\pmb{\Bar{F}} = \Delta\pmb{F}_d (\Delta\pmb{F}_0)_v = [\frac{\text{det}(\Delta\pmb{F}_0)}{\text{det}(\Delta\pmb{F})}]^{1/dof} \Delta\pmb{F} = (\frac{\Delta\Bar{J}}{\Delta{J}})^{1/dof} \Delta\pmb{F},
\end{equation}
where $\Delta{J} = \text{det}(\Delta\pmb{F})$ is the Jacobian of the deformation gradient increment; $\Delta\Bar{J}$ is the averaged Jacobian that has less volumetric constraint compared to $\Delta{J}$. By replacing $\Delta\pmb{F}$ by $\Delta\Bar{\pmb{F}}$ in Eq.~(\ref{eq7}), the volumetric locking can be overcome~\cite{Coombs2018}. 

From the above illustration, the key to the F-bar approach is to determine the averaged Jacobian $\Delta\Bar{J}$ in Eq.~(\ref{eq15}). Departing from the research conducted by Zhao et al.~\cite{Zhao2022}, we calculate the averaged Jacobian $\Delta\Bar{J}$ by mapping the original one $\Delta{J}$ to the related grid nodes and then mapping back to the material points. The mathematical formulation of this mapping and remapping procedure is stated as follows:

(1). The nodal mass-weighted Jacobian can be calculated by mapping the product of mass and $\Delta{J}$ from material points as
\begin{equation}\label{eq16}
	(m\Delta{J})_I = \sum_{p=1}^{{n}_{p}} {N}_{Ip} {m}_{p} \Delta{J},
\end{equation}
where ${m}_{p}$ is the mass of material point $p$; $(m\Delta{J})_I$ is the nodal mass-weighted Jacobian.

(2). Then, the nodal Jacobian can be obtained by
\begin{equation}\label{eq17}
	\Delta{J}_I = \frac{(m\Delta{J})_I}{m_I},
\end{equation}
where ${m_I} = \sum_{p=1}^{{n}_{p}} {N}_{Ip} {m_p}$ is the nodal mass. 

(3). Finally, the nodal Jacobian is mapping back to the material points via
\begin{equation}\label{eq18}
	\Delta\Bar{J} = \sum_{I=1}^{{n}_{n}} {N}_{Ip} \Delta{J}_I.
\end{equation}

Although the above F-bar method can overcome the volumetric locking, the spurious stress oscillation cannot be entirely cured in MPM. In this research, in addition to modifying $\Delta\pmb{F}$, we apply the same mapping remapping technique to the last converged left elastic Cauchy-Green strain ${\pmb{b}}^e(t)$. As a result, the specious stress oscillation can be cured by averaging both $\Delta\pmb{F}$ and ${\pmb{b}}^e(t)$. We name this method \emph{the modified F-bar method.} In this method, Eq.~(\ref{eq7}) is replaced by
\begin{equation}\label{eq19}
	\pmb{b}^{e}(t+\text{d}t) = \Delta\pmb{\Bar{F}} \Bar{\pmb{b}}^e (t) (\Delta\pmb{\Bar{F}})^T,
\end{equation}
where $\Bar{\pmb{b}}^e$ is the averaged left elastic Cauchy-Green strain, which is obtained by
\begin{equation}\label{eq20}
	\Bar{\pmb{b}}^e = \sum_{I=1}^{{n}_{n}} {N}_{Ip} (\frac{1}{m_I}\sum_{p=1}^{{n}_{p}} {N}_{Ip} {m}_{p} \pmb{b}^e).
\end{equation}
Note that Eq.~(\ref{eq20}) represents the same mass-weighted mapping and remapping procedure as Eq.~(\ref{eq16}) to (\ref{eq18}) but in a compact format. {Applying Eq.~(\ref{eq20}) only (i.e. without utilising the F-bar) can also mitigate some form of volumetric locking. However, similarly to solely applying the F-bar method, the stress oscillation cannot be controlled at a satisfactory level.}

This mass-weighted mapping and remapping technique has been applied in the post-processing stage to smooth the contour map of stresses in MPM~\cite{Andersen2013, Dunatunga2015}. {The rationale behind using this mapping and re-mapping technique is that the velocity field of the original MPM is stable despite the highly oscillated stress field because the velocity of a material point has been mass projected to the grids and mapped back to the material point in an MPM computational circle.} However, applying this technique only in the post-processing stage cannot prevent the unphysical oscillation of the stress field. 

The idea of our modified F-bar method is similar to the method proposed by Zhao et al.~\cite{Zhao2022}, but these are different methods because (a) we use a different incremental version of the F-bar method; (b) we apply the mapping and remapping technique to the term $\Delta{J}$ instead of $J$; and (c) we also smooth the left elastic Cauchy-Green strain. These modifications improve the algorithm’s performance under extremely large deformation situations, which will be shown in the later numerical examples. A reader may refer to Zhao et al.~\cite{Zhao2022} for a detailed explanation of the original method.

\subsection{Proposed implicit B-spline MPM algorithm}
\label{subsec:24}
This research follows the implicit MPM algorithm proposed by Guilkey and Weiss~\cite{Guilkey2003}, based on the unconditionally stable Newmark integration. This implicit algorithm was designed for the original MPM. We tailored this implicit MPM algorithm to fit our proposed methods provided in the previous sections. {In the original implicit MPM algorithm~\cite{Guilkey2003}, the Newton–Raphson method is adopted to solve the non-linear equations. However, in this research, the quasi-Newton method proposed by Broyden~\cite{Broyden1965} is used. It allows neglecting the linearisation of the stiffness contribution due to the proposed F-bar method. Therefore, quadratic convergence cannot be achieved. 
Instead, the employment of the aforementioned quasi-Newton method results in a larger number of iterations but less computational time for each iteration.} The following algorithm describes a loading step from time $t$ to $t+\text{d}t$ for the proposed implicit locking-free BSMPM.

\begin{enumerate}
	\item[(1)] \emph{Shape functions and topology}
	
	Given the coordinates of material points $\pmb{x}_p (t)$ and the global knot vectors, construct the B-spline shape functions $N_{Ip}$ and derivatives with respect to global coordinates at the start of a loadstep $\nabla N_{Ip}$ using the algorithm mentioned in Appendix \ref{ap:1}.

	\item[(2)] \emph{Map information from material point to tensor product grid nodes}
	
	(a). The nodal mass $m_I$ can be obtained by 
	
	$m_I = \sum_{p=1}^{{n}_{p}} {N}_{Ip} {m}_{p}$.
	
	(b). The nodal velocity $\pmb{v}_I(t)$ can be obtained by 
	
	$\pmb{v}_I(t) = \frac{1}{m_I}\sum_{p=1}^{{n}_{p}} {N}_{Ip} \pmb{v}_p(t) {m}_{p}$.
	
	(c). The nodal acceleration $\pmb{a}_I(t)$ can be obtained by 
	
	$\pmb{a}_I(t) = \frac{1}{m_I}\sum_{p=1}^{{n}_{p}} {N}_{Ip} \pmb{a}_p(t) {m}_{p}$, 
	
	where $\pmb{a}_p(t)$ is the material point acceleration.
	
	(d). The nodal external force $\pmb{Fext}_I (t+\text{d}t)$ can be obtained by 
	
	$\pmb{Fext}_I = \sum_{p=1}^{{n}_{p}} {N}_{Ip} \pmb{Fext}_p(t+\text{d}t)$, 
	
	where $\pmb{Fext}_p(t+\text{d}t)$ is the material point external force.
	
	\item[(3)] {\emph{Initilisation of the quasi-Newton iteration}}
	
	{(a). Set the number of iterations $k=0$.}
		
	{(b). Use the last converged nodal displacement as the initial guess: $\pmb{u}_I^{k} (t+\text{d}t) = \pmb{u}_I^{k} (t)$.}
	
	\item[(4)] \emph{Update the permanent information at each material point}
	
	(a). Given $\pmb{u}_I^{k} (t+\text{d}t)$ calculate the incremental deformation gradient by 
	
	$\Delta\pmb{F}_p^{k} (t+\text{d}t) = \pmb{I} + \sum_{I=1}^{{n}_{n}} \pmb{u}_I^{k} (t+\text{d}t) \otimes \nabla N_{Ip}$,

	(b). Given $\Delta\pmb{F}_p^{k} (t+\text{d}t)$ and last converged deformation gradient from the last time step $\pmb{F}_p (t)$, update the deformation gradient using 
	
	$\pmb{F}_p^{k} (t+\text{d}t) = \Delta\pmb{F}_p^{k} (t+\text{d}t) \pmb{F}_p (t)$.
	
	(c). Compute the determinants of the deformation gradient and incremental deformation gradient by 
	
	${J}_p^{k} (t+\text{d}t) = \text{det}[\pmb{F}_p^{k} (t+\text{d}t)]$ \& $\Delta{J}_p^{k} (t+\text{d}t) = \text{det}[\Delta\pmb{F}_p^{k} (t+\text{d}t)]$.
	
	(d). Update the volume of material point:
	
	${V}_p(t+\text{d}t) = \Delta{J}_p^{k}(t+\text{d}t) {V}_p (t)$
	
	(e). Compute the averaged Jacobian increment: 
	
	$\Delta\Bar{J}_p^{k}(t+\text{d}t) = \sum_{I=1}^{{n}_{n}} {N}_{Ip} [\frac{1}{m_I}\sum_{p=1}^{{n}_{p}} {N}_{Ip} {m}_{p} \Delta{J}_p^{k}(t+\text{d}t)]$.
	
	(f). Compute the incremental F-bar by
	
	$\Delta\pmb{\Bar{F}}_p^{k}(t+\text{d}t) = [\frac{\Delta\Bar{J}_p^{k}(t+\text{d}t)}{\Delta{J}_p^{k}(t+\text{d}t)}]^{1/dof}\Delta\pmb{F}_p^{k} (t+\text{d}t)$.
	
	(g). Given the incremental F-bar $\Delta\pmb{\Bar{F}}_p^{k}(t+\text{d}t)$ and the last converged averaged left elastic Cauchy-Green strain $\Bar{\pmb{b}}_p^e (t)$, update the left elastic Cauchy-Green strain $\pmb{b}_{p}^{e,k}(t+\text{d}t)$ by
	
	$\pmb{b}_{p}^{e,k}(t+\text{d}t) = \Delta\pmb{\Bar{F}}_p^{k}(t+\text{d}t) \Bar{\pmb{b}}_p^e (t) \Delta\pmb{\Bar{F}}_p^{k}(t+\text{d}t)^T$
	
	(h). Evaluate the elastic logarithmic strain $\pmb{\epsilon}_p^{e,k} (t+\text{d}t)$ by
	
	$\pmb{\epsilon}_p^{e,k} (t+\text{d}t) = \frac{1}{2}\text{ln}\: \pmb{b}_{p}^{e,k}(t+\text{d}t) $
	
	(i). Utilising the return-mapping algorithm, update the Kirchhoff stress $\pmb{\tau}_p^{k}(t+\text{d}t)$ and the left elastic Cauchy-Green strain $\pmb{b}_p^{e,k} (t+\text{d}t)$ based on $\pmb{\epsilon}_p^{e,k} (t+\text{d}t)$. Then, Restore the Cauchy stress $\pmb{\sigma}_p^{k}(t+\text{d}t)$ by
	
	$\pmb{\sigma}_p^{k}(t+\text{d}t) = \frac{\pmb{\tau}_p^{k}(t+\text{d}t)}{J_p^{k}(t+\text{d}t)}$

	\item[(5)] {\emph{Obtain global out-of-balance force matrix}}
	
	{(a). The global nodal internal force $\pmb{Fint}^k_I(t+\text{d}t)$ is obtained by}
	
	{$\pmb{Fint}^k_I(t+\text{d}t) = \sum_{p=1}^{{n}_{p}} \pmb{G}^T \pmb{\sigma}_p^{k}(t+\text{d}t) V_p(t+\text{d}t)$, }
	
	{where $\pmb{G}=\nabla_xN_{Ip} =  \nabla N_{Ip}[\Delta\pmb{F}_p^{k} (t+\text{d}t)]^{-1}$ is the strain-displacement matrix. }
 
        {(b). The out-of-balance force $\pmb{Foob}^{k+1}_I(t+\text{d}t)$ is calculated by}

        {$\pmb{Foob}^{k+1}_I(t+\text{d}t) = \pmb{Fext}_I(t+\text{d}t) - \pmb{Fint}^k_I(t+\text{d}t) - \pmb{M}_I[\frac{4}{\text{d}t^2}\pmb{u}^{k}_I(t+\text{d}t) - \frac{4}{\text{d}t}\pmb{v}_I(t) - \pmb{a}_I(t)]$}

        {where $\pmb{M}_I$ is the nodal mass matrix. The lumped nodal mass $\pmb{M}_I = \text{diag}(m_I)$ is typically used, which is a diagonal matrix.}

        \item[(6)] {\emph{Obtain global nodal stiffness matrix}}

        {if k = 0, get the initial guess of the global stiffness matrix by}

        {$\pmb{K}^k_I(t+\text{d}t) = \frac{4}{\text{d}t^2}\pmb{M}_I + \sum_{p=1}^{{n}_{p}} \pmb{G}^T \pmb{D}_p^{e,k}(t+\text{d}t) \, \pmb{G} \, V_p(t+\text{d}t)$}

        {where, $\pmb{D}_p^{e,k}(t+\text{d}t)$ is the elastic tangent modulus of a material point.}

        {if k > 0, estimate the global stiffness matrix by}

        {$\pmb{K}^k_I(t+\text{d}t) = \pmb{K}^k_I(t) - \frac{\pmb{Foob}^{k+1}_I(t+\text{d}t) {\delta\pmb{u}^{k}_I}^T}{{\delta\pmb{u}^{k}_I}^T \delta\pmb{u}^{k}_I} $}

        {where, $\delta\pmb{u}^{k}_I$ is the nodal displacement increment from the last quasi-Newton iteration.}

        \item[(7)] \emph{Solve the linear equation}

        At {quasi-Newton} iteration number $k$, the nodal displacement increment $\delta\pmb{u}^{k+1}_I$ can be obtained by solving the following linear equation with Dirichlet boundary conditions

        $\pmb{K}^k_I(t+\text{d}t) \cdot \delta\pmb{u}^{k+1}_I =  \pmb{Foob}^{k+1}_I(t+\text{d}t)$

        \item[(8)] \emph{Update grid nodal displacement:}

        $\pmb{u}^{k+1}_I(t+\text{d}t) = \pmb{u}^{k}_I(t+\text{d}t) + \delta\pmb{u}^{k+1}_I$

        \item[(9)] \emph{Check convergency}

        If $f^{k+1} = \frac{\|\delta\pmb{u}^{k+1}_I\pmb{Foob}^{k+1}_I\|}{\|\delta\pmb{u}^{k=1}_I\pmb{Foob}^{k=1}_I\|} < tol$

        go to step (10). Otherwise, set $k=k+1$ and go to step (4).

        \item[(10)] \emph{Update the grid kinematics}

        (a). The nodal acceleration $\pmb{a}_I(t+\text{d}t)$ can be updated by

        $\pmb{a}_I(t+\text{d}t) = \frac{4}{\text{d}t^2}\pmb{u}^{k+1}_I(t+\text{d}t) - \frac{4}{\text{d}t}\pmb{v}_I(t) - \pmb{a}_I(t)$

        (b). The nodal velocity $\pmb{v}_I (t+\text{d}t)$ can be updated by

        $\pmb{v}_I(t+\text{d}t) = \frac{2}{\text{d}t}\pmb{u}^{k+1}_I(t+\text{d}t) - \pmb{v}_I(t)$

        \item[(11)] \emph{Save the converged information into the material points}

        (a). The spatial coordinates of material points $\pmb{x}_p(t+\text{d}t)$ can be updated by

        $\pmb{x}_p(t+\text{d}t) = \pmb{x}_p(t) + \sum_{I=1}^{{n}_{n}} {N}_{Ip} \, \pmb{u}^{k+1}_I(t+\text{d}t)$

        (b). The displacement of material points $\pmb{u}_p(t+\text{d}t)$ can be updated by

        $\pmb{u}_p(t+\text{d}t) = \pmb{u}_p(t) + \sum_{I=1}^{{n}_{n}} {N}_{Ip} \, \pmb{u}^{k+1}_I(t+\text{d}t)$

        (c). The acceleration of material points $\pmb{a}_p(t+\text{d}t)$ can be updated by

        $\pmb{a}_p(t+\text{d}t) = \sum_{I=1}^{{n}_{n}} {N}_{Ip} \, \pmb{a}_I(t+\text{d}t)$

        (d). The velocity  of material points $\pmb{v}_p(t+\text{d}t)$ can be updated by

        $\pmb{v}_p(t+\text{d}t) = \pmb{v}_p(t) + \frac{\pmb{a}_p(t)+\pmb{a}_p(t+\text{d}t)}{2}\text{d}t$

        (e). Given the updated left elastic Cauchy-Green strain $\pmb{b}_p^{e,k} (t+\text{d}t)$, compute the averaged left elastic Cauchy-Green strain $\Bar{\pmb{b}}_p^{e} (t+\text{d}t)$ by
	
	$\Bar{\pmb{b}}_p^e (t+\text{d}t) = \sum_{I=1}^{{n}_{n}} {N}_{Ip} [\frac{1}{m_I}\sum_{p=1}^{{n}_{p}} {N}_{Ip} {m}_{p} \pmb{b}_p^{e,k} (t+\text{d}t)]$
        
        (f). Finally, save all the converged permanent variables like $\pmb{F}_p(t+\text{d}t)$ and $V_p(t+\text{d}t)$, then go to the next time step.
\end{enumerate}

\section{Numerical examples}
\label{sec:3}

Both dynamic and {quasi-static} large deformation numerical examples are included in this research to illustrate the performance of the new method proposed in this paper. We set the mass matrix to zero for the {quasi-static} analysis, and hence the dynamic implicit algorithm automatically becomes {quasi-static}~\cite{Guilkey2003}. As mentioned before, the Mohr-Coulomb constitutive model has been applied in this research. The third numerical example is based on the Tresca yield criterion, which is a special case of Mohr-Coulomb with zero friction and dilation angles. Also, the soil is assumed as weightless for the third example. The material input parameters for these numerical models are summarised in table~\ref{tab:1}. {The dilation angle for these three examples is zero}. The simulations are conducted using an in-house developed MATLAB code. The VTK files are generated during the simulations, and the post-processing of these VTK files is performed using the software ParaView.

\begin{table}[!t]
\centering
\caption{ Material properties for three numerical examples.} \label{tab:1}
	\vspace*{0.2cm}
	\begin{tabular}{c|ccccc}
	\hline\noalign{\smallskip}
	Numerical example & Unit weight & Young’s modulus & Poisson’s ratio & Friction angle & Cohesion \\
        Section & $\gamma[kN/m^3]$ & $E[kPa]$ & $\nu$ & $\phi[^\circ]$ & $c[kPa]$ \\
		\noalign{\smallskip}\hline\noalign{\smallskip}
	\ref{subsec:31} & 20.4 & 5840    & 0.3   & 21.9 & 0  \\
	\ref{subsec:32} & 20   & 100,000 & 0.3   & 20   & 10 \\
        \ref{subsec:33} & 0    & 100     & 0.495 & 0 & 1  \\
\hline\noalign{\smallskip}
	\end{tabular}
\end{table}

\subsection{Granular column collapse}
\label{subsec:31}

The granular column collapse experiment was initially designed to study natural catastrophe granular flow propagation~\cite{Lagree2011}. This research uses a granular column collapse experiment documented by Nguyen et al.~\cite{Nguyen2017} in order to validate the proposed numerical model. Nguyen et al.~\cite{Nguyen2017} use a Smoothed Particle Hydrodynamics (SPH) method with Mohr-Coulomb constitutive model to validate this experiment, where aluminium rods are used to represent the granular soil. The material properties of these aluminium rods were calibrated using a series of experimental tests~\cite{Nguyen2017}, summarised in table~\ref{tab:1}. 

\begin{figure}
\begin{center}
\includegraphics[width=10cm]{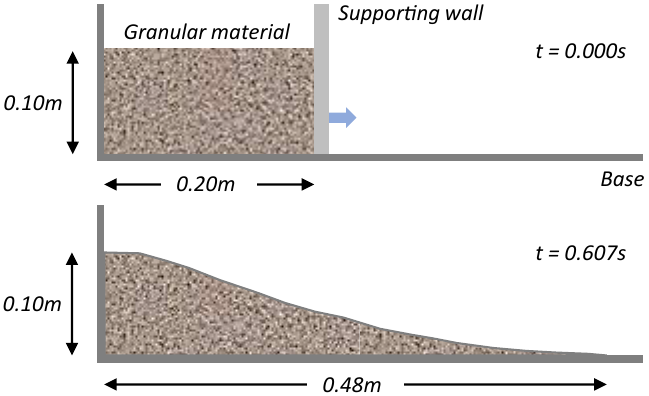}
\caption{Graphical illustration of granular column collapse experiment.}
\label{fig2}
\end{center}
\end{figure}

Fig.\ref{fig2} graphically illustrates this granular column collapse experiment. Initially, the granular materials were set to be rectangular, with 0.1 m in height and 0.2 m in width~\cite{Nguyen2017}. Then, the granular flow was initialised by suddenly removing the supporting wall. The evolution of the granular flow was recorded using a high-speed camera. Finally, the granular flow became stationary at 0.607 seconds after removing the supporting wall. 

The aluminium rods have neglectable movement in the z-direction (in page direction), which allows us to model it under a 2-D plane strain condition. In the 2-D numerical model, the supporting wall that blocks the aluminium bars shown in Fig.\ref{fig2} has not been modelled. The gravity acceleration is applied at the initial time step since the wall was removed quickly during the experiment~\cite{Nguyen2017}. The width and height of the background mesh (parametric grid) are set to be 550 and 102.5 mm, respectively. Both the width and height for each parametric grid cell (patch) are 2.5 mm. Thirty-six ($6^2$) material points are set in each grid cell to minimise the quadrature error. As a result, the aluminium rods have been discretised into 115,200 material points. The tensor product grid is formed from the parametric grid. The fixed and roller boundary conditions are applied on the bottom and sides of the grid, respectively. 

Viscous damping was applied by Nguyen et al.~\cite{Nguyen2017} to stabilise their SPH-based model. In this research, no numerical damping is applied because the energy loss during the rotation of the aluminium rods is negligible~\cite{Solowski2015}. Finally, the time increment is set as 0.00025 seconds, resulting in 2,428 steps in total. The tolerance for the Newton-Raphson iteration is set to be $10^{-8}$.

\begin{figure}
\begin{center}
\includegraphics[width=\textwidth]{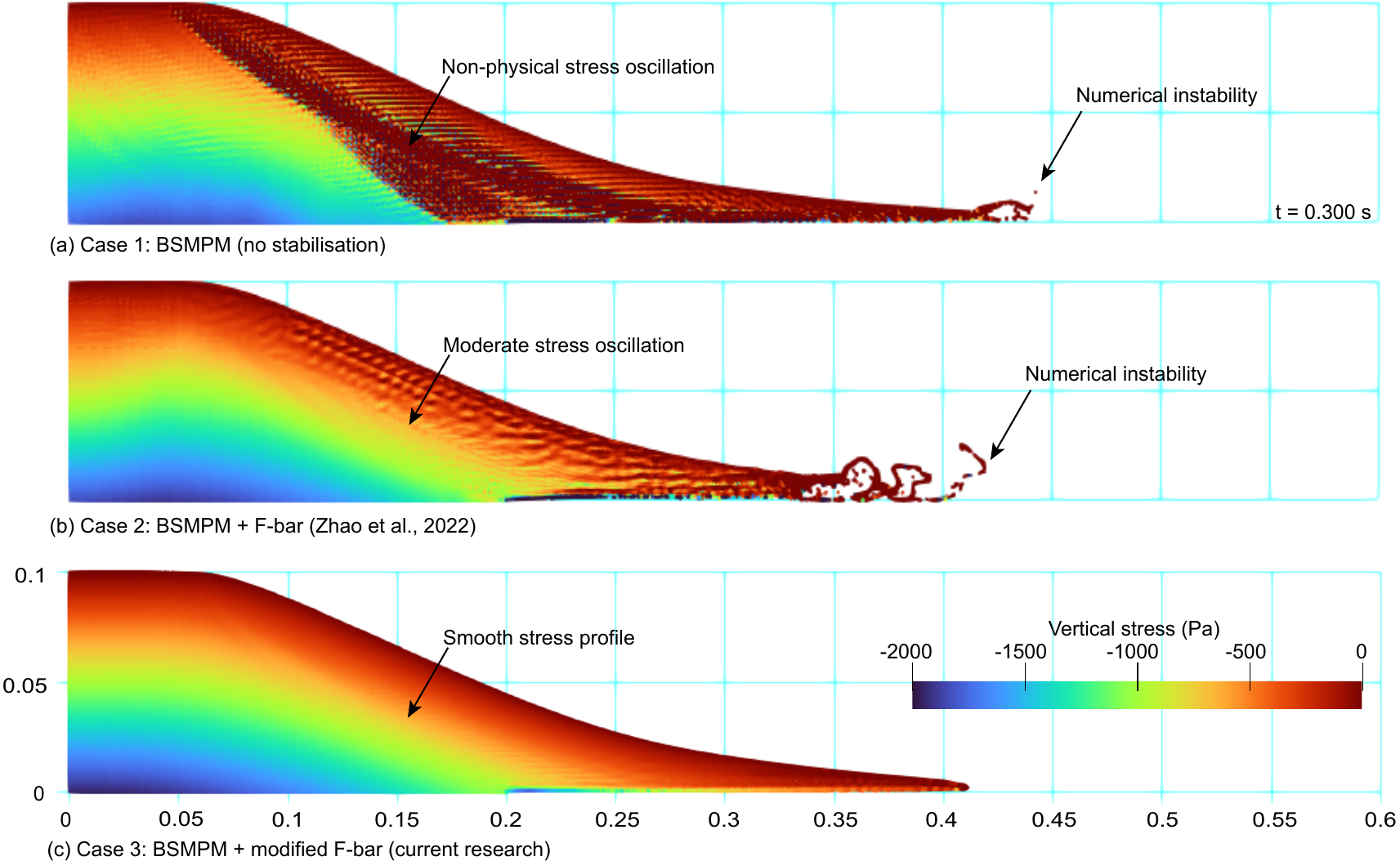}
\caption{Simulations of granular column collapse by different methods: (a) BSMPM with no stabilisation; (b) BSMPM with F-bar method; (c) BSMPM with modified F-bar method.}
\label{fig3}
\end{center}
\end{figure}

In this numerical example, three scenarios have been examined: the BSMPM with no stabilisation method, the BSMPM with F-bar method (Zhao et al. ~\cite{Zhao2022}) and BSMPM with the newly proposed modified F-bar method. The first two methods experience numerical instability after the granular flow reaches a relatively long runout distance (around 0.37 m). The second one, the BSMPM with the F-bar method (Zhao et al. ~\cite{Zhao2022}), cannot have a converged result after 0.351 seconds due to the significant numerical instability. The first simulation, the BSMPM with no stabilisation, has a converged result, but the unstable granular flow clashes with the right boundary (i.e. 0.55 m) at about 0.42 seconds. These findings are illustrated in Fig.\ref{fig3}, which shows the vertical stress contour maps at 0.300 seconds given by these three methods. Fig.\ref{fig3}(a) shows that the BSMPM with no stabilisation experiences severe stress oscillation due to the kinematic locking. The over-stiff behaviour is not observed in the BSMPM for this granular column collapse problem. However, the non-physical stress oscillation causes kinematic instability and over-prediction of the runout distance. Fig.\ref{fig3}(b) shows that the F-bar method (Zhao et al. ~\cite{Zhao2022}) significantly reduces stress oscillation. However, we can still observe moderate stress oscillation and amplified kinematic instability, as shown in Fig.\ref{fig3}(b). Fig.\ref{fig3}(c) shows that only minimal stress oscillation is observed in our modified F-bar method, and the stress profile is smooth. The remaining stress oscillation can be prevented by reducing the time step, resulting in losing computational efficiency. This research uses the relatively large time step by trading off a limited amount of accuracy. 

\begin{figure}
	\begin{center}
		\includegraphics[width=10cm]{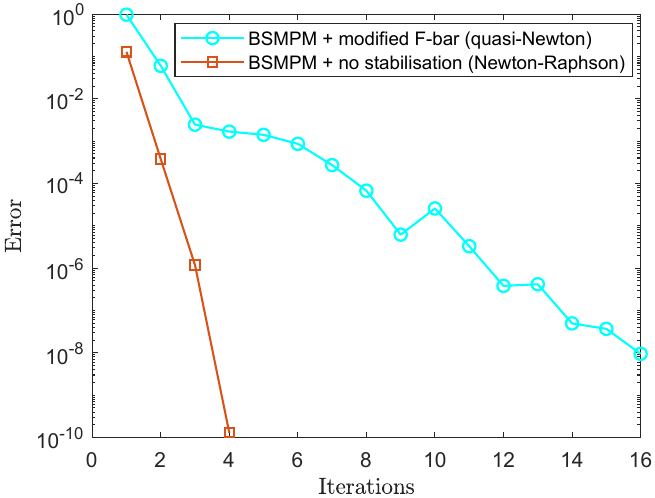}
		\caption{Granular column collapse, load step 10. Comparison of quasi-Newton and Newton-Raphson methods }
		\label{convrate}
	\end{center}
\end{figure}

{ Fig.\ref{convrate} shows the comparison of quasi-Newton and Newton-Raphson methods. A linear convergence can be observed in this quasi-Newton method, while the convergence rate for the Newton-Raphson method is quadratic. However, the computational effort for each iteration is very small for this quasi-Newton method. The overall efficiency of the quasi-Newton method can be higher than the Newton-Raphson method, especially for a highly non-linear large-scale problem. }

\begin{figure}
\begin{center}
\includegraphics[width=10cm]{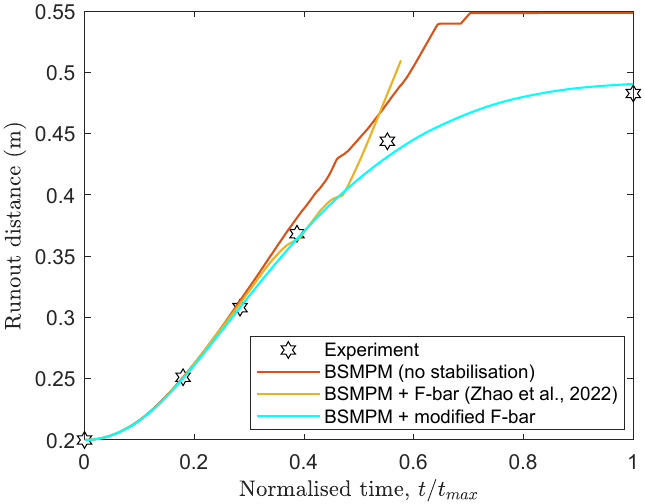}
\caption{Evolutions of the runout distance.}
\label{fig4}
\end{center}
\end{figure}

\begin{figure}
\begin{center}
\includegraphics[width=\textwidth]{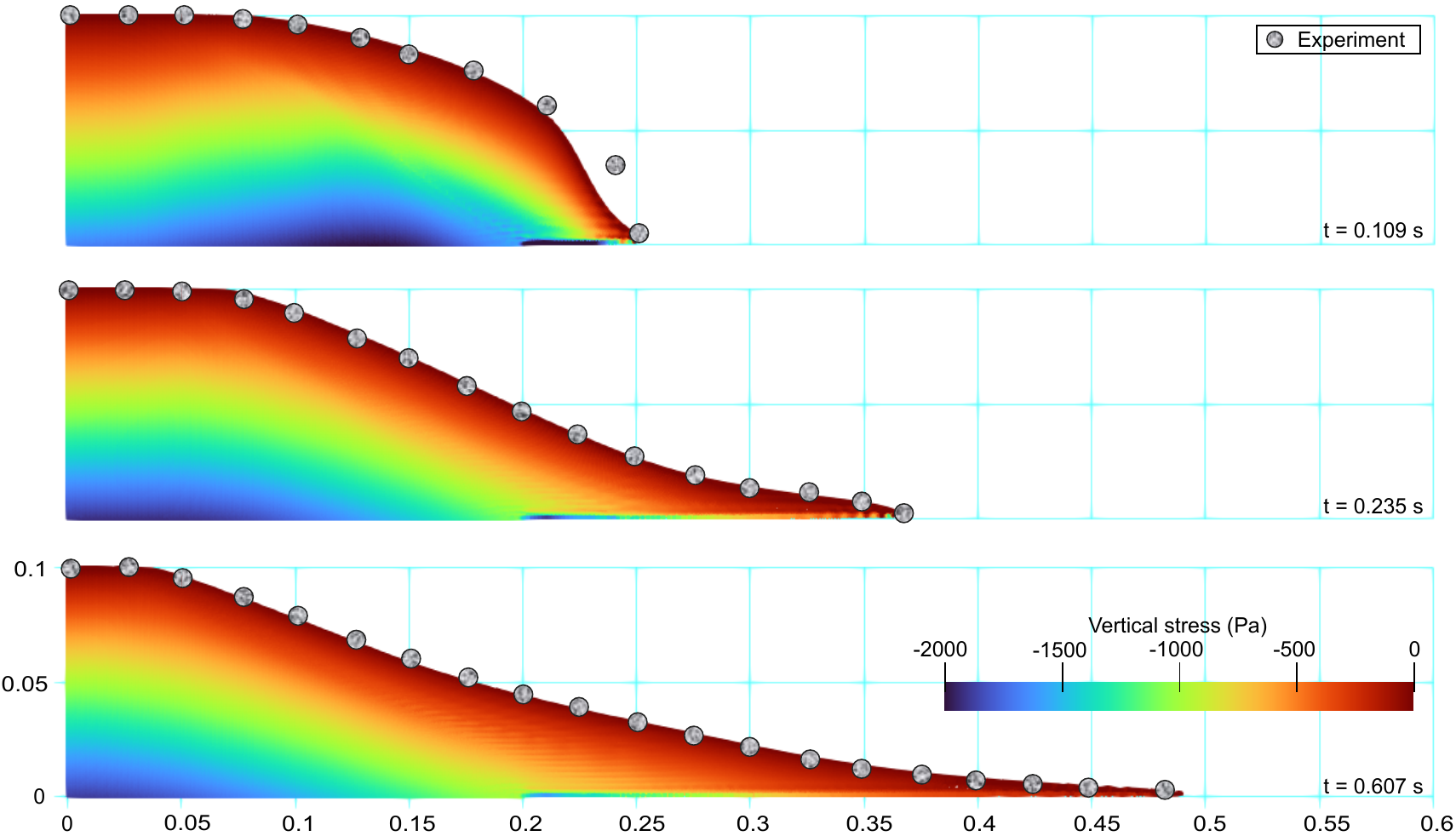}
\caption{Comparison between the BSMPM with modified F-bar method with the experiment.}
\label{fig5}
\end{center}
\end{figure}

Fig.\ref{fig4} compares the evolutions of the runout distance. From Fig.\ref{fig4}, only our approach has excellent agreement with the experiment. Fig.\ref{fig5} compares the free surface of granular flow given by our approach and the experiment at 0.109, 0.235 and 0.607 seconds, respectively. The grey dots in Fig.\ref{fig5} are the digitised experimental results, representing the free surface of granular flow. Fig.\ref{fig5} shows that the proposed method has a remarkable agreement with the experiment during the entire propagation process of the granular flow. 

{There is discontinuity of the stress field on the base of the simulation between 0.2 and 0.25 in Fig.\ref{fig3} and Fig.\ref{fig5}. This discontinuity is caused by the material point moving to the fixed boundary condition on the bottom surface. Modelling a frictional contact can solve this issue and result in a more realistic result. However, in this research, the fixed boundary condition is used to represent the rough surface for simplicity.}

\subsection{Slope stability}
\label{subsec:32}

The second numerical example is a static slope stability problem introduced by Griffiths and Lane~\cite{Griffiths1999}, who used the shear strength reduction technique to determine the factor of safety (FOS) of the slope. In the strength reduction method, the FOS of a slope can be represented by the ratio between the actual soil strength and the reduced one. Therefore, the factored shear strength parameter (i.e. cohesion and friction angle) $c_f$ and $\phi_f$ can be determined by~\cite{Griffiths1999}: $c_f = c/SRF$ and $\phi_f = \text{arctan}(\frac{\text{tan}\phi}{SRF})$, where SRF is the shear strength reduction factor (SRF). The critical value of SRF represents the FOS of the slope.

\begin{figure}
\begin{center}
\includegraphics[width=10cm]{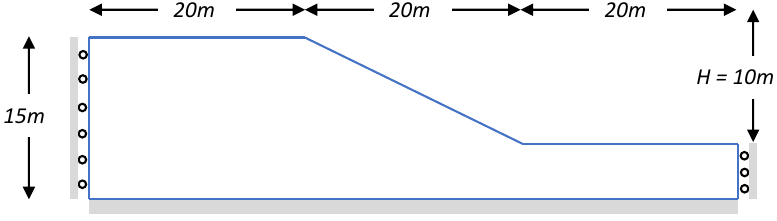}
\caption{Geometry and boundary conditions of the slope example.}
\label{fig6}
\end{center}
\end{figure}

The geometry of this slope is shown in Fig.\ref{fig6}, and its material properties are summarised in table~\ref{tab:1}. As shown in Fig.\ref{fig6}, the height of the slope $H$ is 10 meters, and the inclination angle is 26.57 degrees. We keep the dimensionless soil property $c/(\gamma H)=0.05$ same as the one used by Griffiths and Lane~\cite{Griffiths1999}. The fixed and roller boundary conditions are applied on the bottom and sides of the slope, respectively. 

\begin{figure}
\begin{center}
\includegraphics[width=10cm]{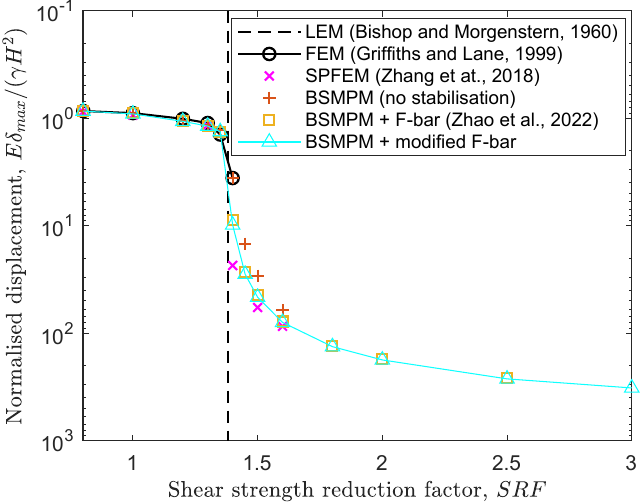}
\caption{Normalised maximum displacement with different SRFs by different methods.}
\label{fig7}
\end{center}
\end{figure}

Griffiths and Lane~\cite{Griffiths1999} analysed this slope using small deformation FEM up to SRF = 1.4. They plotted the SRF versus $E\delta_{max}/(\gamma H^2)$ (a dimensionless displacement), as shown in Fig.\ref{fig7}. $\delta_{max}$ is the maximum displacement in this slope. However, due to the highly distorted mesh, the FEM yields an unconverged solution at SRF = 1.4, and the displacement at SRF = 1.4 in Fig.\ref{fig7} results from the unconverged 1,000 iterations~\cite{Griffiths1999}. In FEM, the failure of the slope is defined as the stage that the algorithm cannot find a converged solution within a user-defined number of iterations, which indicates the FOS of this slope is between 1.35 and 1.4~\cite{Griffiths1999}. This FOS is consistent with the FOS = 1.38 determined by Bishop and Morgenstern~\cite{Bishop1960} using the limit equilibrium method (LEM). 

In this research, we follow the procedure that Griffiths and Lane~\cite{Griffiths1999} proposed to determine the dimensionless displacement $E\delta_{max}/(\gamma H^2)$ of this slope up to SRF = 3.0. The size of the parametric grid cell (patch) is 1 $m^2$, which is similar to the mesh size used in Griffiths and Lane~\cite{Griffiths1999}. Same as the previous example, thirty-six ($6^2$) material points are used in each grid cell. As a result, the slope is discretised by 21,540 material points. Each SRF represents an independent analysis. For the BSMPM with F-bar methods, the gravity load gradually increases from zeros to the maximum value in 10 steps for SRF < 1.38 and 100 steps for SRF > 1.38, respectively. For the BSMPM without stabilisation, 500 steps of loading are required after the failure (SRF > 1.38) to ensure numerical stability. In this numerical example, no kinematic instability has been observed. However, the BSMPM with no stabilisation cannot have a converged solution after SRF = 1.6, and the BSMPM with the F-bar method (Zhao et al.~\cite{Zhao2022} cannot have a converged result after SRF = 2.5.

\begin{figure}
\begin{center}
\includegraphics[width=\textwidth]{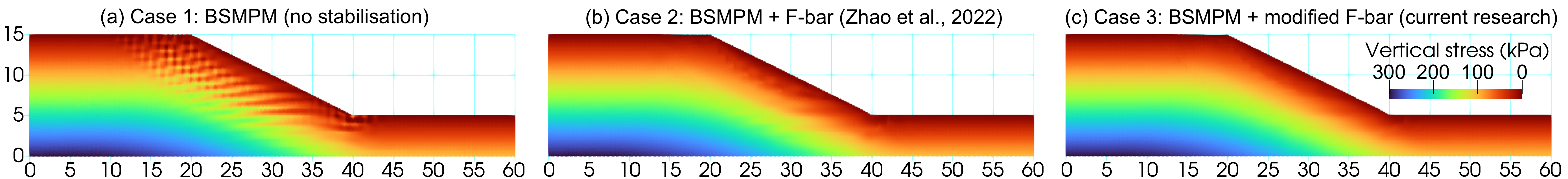}
\caption{Simulations of the slope stability at SRF = 1.4 by different methods: (a) BSMPM with no stabilisation; (b) BSMPM with F-bar method; (c) BSMPM with modified F-bar method.}
\label{fig8}
\end{center}
\end{figure}

Before the slope failure (SRF < 1.38), the results of three BSMPMs closely match the FEM research of Griffiths and Lane~\cite{Griffiths1999}, as shown in Fig.\ref{fig7}. After the failure, the BSMPM without stabilisation methods generates slightly over-stiff results due to the volumetric locking. We can observe significant stress oscillation in the contour map generated by BSMPM, as shown in Fig.\ref{fig8}(a). Similar to the previous numerical example, the F-bar method Zhao et al.~\cite{Zhao2022} proposed significantly reduces the volumetric locking and stress oscillation. In contrast, the modified F-bar method completely overcomes the volumetric locking and stress oscillation. Fig.\ref{fig7} shows that the BSMPM (no stabilisation) has a similar result as the FEM when SRF = 1.4. However, both results are inaccurate because of either the volumetric locking instability in BSMPM (Fig.\ref{fig8}(a)) or the highly distorted FEM mesh~\cite{Griffiths1999}. Because of the volumetric locking instability, the BSMPM cannot have a converged solution after SRF = 1.6, even further reducing the time step. Although the maximum displacements given by the two F-bar methods are similar, the newly proposed modified F-bar method results in a slightly larger displacement for this numerical example, with a completely smooth stress solution, as illustrated in Fig.\ref{fig8}(c). 

Zhang et al.~\cite{Zhang2018} investigated the same slope problem using the smooth particle finite element method (SPFEM), a large deformation finite element-based method. The results of BSMPM with F-bar methods are stiffer than Zhang et al.~\cite{Zhang2018} because the small strain constitutive model is used in their research, overestimating the displacement in the large strain range. This phenomenon is also reported by Borja and Tamagnini~\cite{Borja1998} and Coombs and Crouch~\cite{Coombs2011b}.

\begin{figure}
\begin{center}
\includegraphics[width=\textwidth]{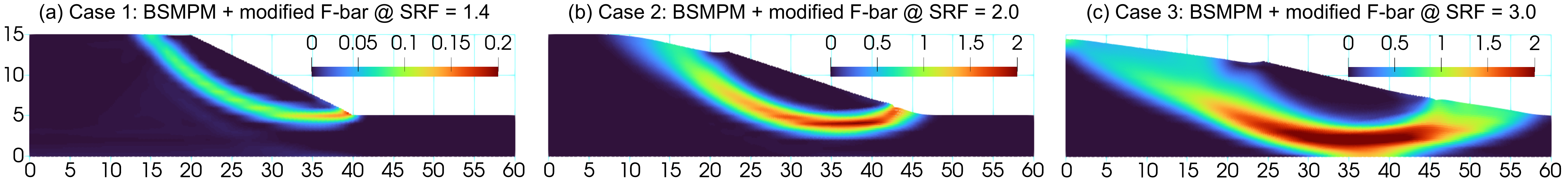}
\caption{Evolution of failure surface: the accumulated plastic strain at different SRFs: (a)SRF=1.4; (b) SRF=2.0; (c) SRF=3.0.}
\label{fig9}
\end{center}
\end{figure}

Fig.\ref{fig9} shows the evolution of the failure surface at the post-failure stages (SRF > 1.38). The value plotted in Fig.\ref{fig9} is the accumulated plastic strain generated by the BSMPM with the modified F-bar method. The accumulated plastic strain is an internal variable that describes the development of shear plasticity and controls the Mohr-Coulomb model’s hardening behaviour. However, we assume no hardening behaviour (perfectly plastic) in these numerical examples (a detailed explanation of this accumulated plastic strain and the hardening behaviour is reported by de Souza Neto et al.~\cite{deSouza2008}). As shown in Fig.\ref{fig9}, the failure surface develops in width as the increasing SRF. Finally, at SRF = 3.0, the failure surface reaches the bottom of the slope foundation.

\subsection{Large penetration of strip footing}
\label{subsec:33}

\begin{figure}
\begin{center}
\includegraphics[width=10cm]{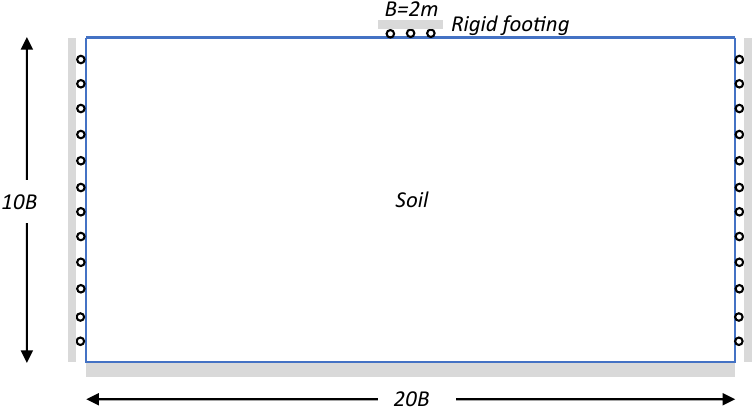}
\caption{Geometry and boundary conditions of the strip footing example.}
\label{fig10}
\end{center}
\end{figure}

The third example is the large penetration of a rigid strip footing on nearly incompressible Tresca soil. We follow the same strip footing example carried out through the Particle Finite Element method (PFEM) by Monforte et al.~\cite{Monforte2017}, who report the mixed stabilised formulation to overcome the volumetric locking. Also, the same large strain constitutive model is used in their research, ensuring the results are comparable with current research. The soil parameters are shown in table~\ref{tab:1}. Fig.\ref{fig10} illustrates the geometry and boundary conditions of this example. We model the rigid strip footing as the prescribed Dirichlet boundary condition in the y-direction (roller) using the penalty method proposed by Cortis et al.~\cite{Cortis2018}. The width of the strip footing is two meters ($B=2 m$), and the final penetration reaches two meters ($z=2 m$) in 160 steps. The internal nodal forces where the Dirichlet boundary condition is prescribed are measured at each time step to represent the pressure underneath the footing. Due to the symmetry, only the right half of the problem has been modelled. The parametric grid cell (patch) size is 0.01 $m^2$, and thirty-six ($6^2$) material points are used in each grid cell resulting in a total of 360,000 material points.

\begin{figure}
\begin{center}
\includegraphics[width=10cm]{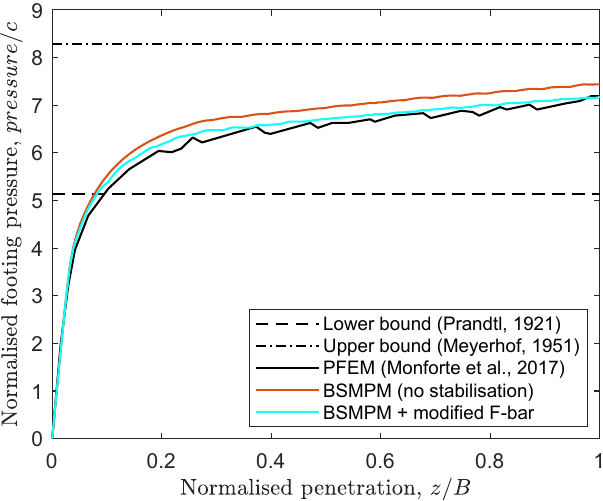}
\caption{Normalised maximum displacement versus normalised pressure underneath the footing by different methods.}
\label{fig11}
\end{center}
\end{figure}

\begin{figure}
\begin{center}
\includegraphics[width=\textwidth]{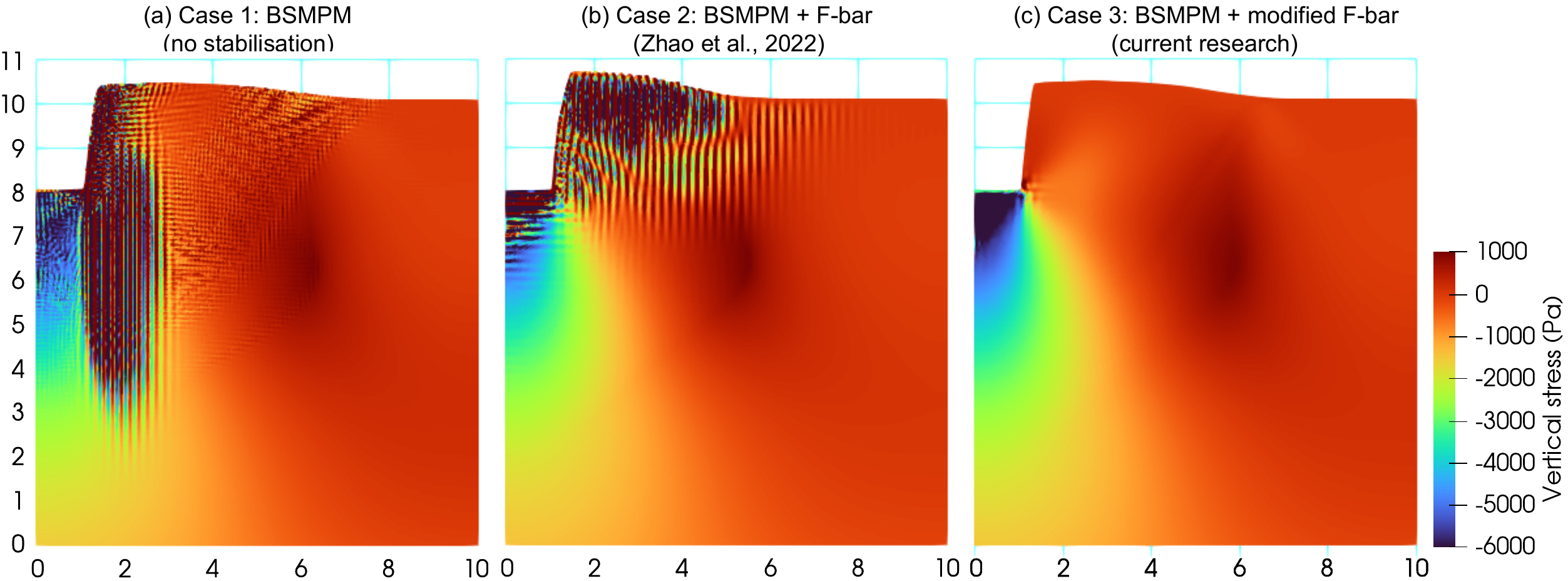}
\caption{Simulations of large penetration of rigid strip footing: (a) BSMPM with no stabilisation; (b) BSMPM with F-bar method; (c) BSMPM with modified F-bar method.}
\label{fig12}
\end{center}
\end{figure}

Fig.\ref{fig11} shows the normalised penetration versus normalised pressure underneath the footing, and Fig.\ref{fig12} shows the vertical stress contours given by the three methods. As shown in Fig.\ref{fig11}, the result given by the proposed method is located between the lower bound, $(\pi+2)c$, and upper bound, $(2\pi+2)c$, analytical solutions derived by Prandtl~\cite{Prandtl1921} and Meyerhof~\cite{Meyerhof1951}, respectively. The modified F-bar (current research) closely matches the PFEM~\cite{Monforte2017} result, as shown in Fig.\ref{fig11}. {We can see a very small amount of fluctuations in the curves because the nodal forces are extracted as the pressure underneath the footing.} On the contrary, the proposed BSMPM with the modified F-bar method shows a less fluctuating result than the PFEM. The BSMPM with no stabilisation has an over-stiff behaviour due to the volumetric locking. The F-bar method (Zhao et al.~\cite{Zhao2022}) has partially solved the volumetric locking issue, but its stress field highly oscillates in the area underneath the footing and the tensile region where the volume expands, as shown in Fig.\ref{fig12}(b). On the contrary, the proposed BSMPM with the modified F-bar method has a very stable stress field {despite some stress oscillation near the corner of the footing}, as shown in Fig.\ref{fig12}(c). {A similar corner issue can be found in the research of Yuan et al.~\cite{Yuan2021}, which may be a common phenomenon for a large deformation strip footing on a nearly incompressible material. We found that slightly reducing Poisson's ratio can help with this issue.}

\section{Conclusions}
\label{sec:4}

In this research, we propose an implicit locking-free BSMPM. The essential ingredients of this BSMPM include: (1) a special BSMPM implementation, (2) a large strain constitutive model and (3) an enhancement of the F-bar method. Thanks to these features, the proposed BSMPM has the following advantages: (a) no cell crossing noise; (b) no volumetric locking; (c) no complex material point searching algorithm; (d) no barrier to switching from the original MPM; (e) reasonably smooth stress profile; and (f) easy application of the volumetric locking strategy.

The proposed numerical model is validated with an experiment and previous numerical studies. These numerical examples disclose the accuracy of the proposed method for extremely large deformation geotechnical problems. The excellent agreement with the granular column collapse experiment indicates that the proposed method can accurately simulate granular flows. The proposed method also shows an accuracy that closely matches the previous FEM slope study in the small deformation range. Furthermore, it can simulate the post-failure stage of the slope and footing problems that involve extensive deformation without the error due to highly distorted FEM mesh. The proposed method can always maintain a smooth and stable stress field, which was a challenge for the previous MPMs. 

The over-stiffening behaviour has been observed in the slope and footing problems because of the volumetric locking. However, in the granular column collapse problem, the stress oscillation due to the volumetric looking leads to a over-estimated runout distance, and no obvious over-stiffening behaviour is shown for the BSMPM with no stabilisation. The volumetric locking instability, found in previous MPMs, results in inaccurate results for all of three numerical examples. Therefore, removing the volumetric locking is essential for a geotechnical application even if the material is not nearly-incompressible.

In future research, it is necessary to extend this method to a soil-water coupled model in order to enlarge the range of geotechnical applications to be solved. Additionally, more advanced constitutive modelling is desired to have more realistic simulations.

\subsection*{Conflict of interest}

The authors declare no potential conflict of interest.





\appendix
\section{Construct the B-spline shape function}
\label{ap:1}

The following scheme illustrates the application of the proposed method for a 2-D problem:

\begin{enumerate}
    \item[(1)] Given the x-coordinates of a material point ${x}_x$ and the global knot vector in x-direction ${\Xi}_{x}$, construct the 1-D B-spline shape functions ${N}_{i,q}(x)$ and derivatives $\frac{\text{d}{N}_{i,q}(x)}{\text{d}x}$ recursively using Eq.~(\ref{eq4}) to ~(\ref{eq6}). Repeat the same process in y-direction to obtain ${N}_{i,q}(y)$ and $\frac{\text{d}{N}_{i,q}(y)}{\text{d}y}$.

    \item[(2)] Convolute shape functions and derivatives in multi-dimensions through the tensor product operation:

    $N_{Ip} = \text{kron}({N}_{i,q}(x), {N}_{i,q}(y))$

    $\partial_x N_{Ip} = \text{kron}(\frac{\text{d}{N}_{i,q}(x)}{\text{d}x}, {N}_{i,q}(y))$

    $\partial_y N_{Ip} = \text{kron}({N}_{i,q}(x), \frac{\text{d}{N}_{i,q}(y)}{\text{d}y})$

    \item[(3)] Obtained the material point to tensor product grid nodes connectivity ($p2N$) by finding the non-zero indexes of the shape functions:
    
    $p2N = \text{find}(N_{Ip} \neq 0)$

    \item[(4)] Chop off the shape functions and derivatives matrices using $p2N$:

    $N_{Ip} = N_{Ip}[p2N]$

    $\nabla N_{Ip} = [\partial_x N_{Ip}[p2N];\; \partial_y N_{Ip}[p2N]]$
    
 \end{enumerate}

In this algorithm, the functions kron() and find() represent the native MATLAB functions or the NumPy functions in Python. The particle searching is achieved by using this find() function. The find(A$\neq$0) function is capable to find the indices of nonzero elements for an input vector A; the kron(A,B) function returns the Kronecker tensor product of two input vectors A and B. The proposed method also works on the non-equally spaced structured mesh by simply utilising non-uniform knot vectors. Moving from 2-D to 3-D is also straightforward by introducing an additional knot vector in the z-direction. The 3-D shape functions can be constructed following the same tensor product procedure.

{The algorithm mentioned above is not efficient because it involves many non-essential zero multiplications (i.e. zero times zero) during the tensor product operations. To accelerate this algorithm, we can chop off (i.e. only keep the non-zero elements) the 1-D shape functions and derivatives before the tensor product operation. Therefore, the matrix size can be reduced significantly during the execution of kron() function. Following the same concept, we can find the non-zero indexes of the 1-D shape function. Then, we can use this 1-D topology information to calculate the multi-dimensions $p2N$ matrix. For 2-D quadratic B-spline MPM, a more efficient algorithm is:}

\begin{enumerate}
    \item[(1)] Given the x-coordinates of a material point ${x}_x$ and the global knot vector in x-direction ${\Xi}_{x}$, construct the 1-D B-spline shape functions ${N}_{i,q}(x)$ and derivatives $\frac{\text{d}{N}_{i,q}(x)}{\text{d}x}$ recursively using Eq.~(\ref{eq4}) to ~(\ref{eq6}). Repeat the same process in y-direction to obtain ${N}_{i,q}(y)$ and $\frac{\text{d}{N}_{i,q}(y)}{\text{d}y}$.

    \item[(2)] Obtain the 1-D topology information $id_x$ in the x-direction and $id_y$ in the y-direction by finding the non-zero indexes of the 1-D shape functions:
    
    $id_x = \text{find}(N_{i,q}(x) \neq 0)$

    $id_y = \text{find}(N_{i,q}(y) \neq 0)$

    \item[(3)] Given $id_x$, $id_y$ and the number of shape functions in y-direction $n_y$, calculate the material point to tensor product grid nodes connectivity ($p2N$):

    $A = (id_x[1]-1)*n_y + id_y$

    $B = A + n_y$

    $C = B + n_y$

    $p2N = [A,\; B,\; C]$

    \item[(4)] Chop off and convolute shape functions and derivatives in multi-dimensions through the tensor product operation:

    $N_{Ip} = \text{kron}({N}_{i,q}(x)[id_x], {N}_{i,q}(y)[id_y])$

    $\partial_x N_{Ip} = \text{kron}(\frac{\text{d}{N}_{i,q}(x)}{\text{d}x}[id_x], {N}_{i,q}(y)[id_y])$

    $\partial_y N_{Ip} = \text{kron}({N}_{i,q}(x)[id_x], \frac{\text{d}{N}_{i,q}(y)}{\text{d}y}[id_y])$

    $\nabla N_{Ip} = [\partial_x N_{Ip};\; \partial_y N_{Ip}]$
    
 \end{enumerate}

 {In this algorithm, $id_x[1]$ means the first element of the $id_x$ vector. For 2-D cubic B-spline MPM, a more efficient algorithm is:}

\begin{enumerate}
    \item[(1)] Given the x-coordinates of a material point ${x}_x$ and the global knot vector in x-direction ${\Xi}_{x}$, construct the 1-D B-spline shape functions ${N}_{i,q}(x)$ and derivatives $\frac{\text{d}{N}_{i,q}(x)}{\text{d}x}$ recursively using Eq.~(\ref{eq4}) to ~(\ref{eq6}). Repeat the same process in y-direction to obtain ${N}_{i,q}(y)$ and $\frac{\text{d}{N}_{i,q}(y)}{\text{d}y}$.

    \item[(2)] Obtain the 1-D topology information $id_x$ in the x-direction and $id_y$ in the y-direction by finding the non-zero indexes of the 1-D shape functions:
    
    $id_x = \text{find}(N_{i,q}(x) \neq 0)$

    $id_y = \text{find}(N_{i,q}(y) \neq 0)$

    \item[(3)] Given $id_x$, $id_y$ and the number of shape functions in y-direction $n_y$, calculate the material point to tensor product grid nodes connectivity ($p2N$):

    $A = (id_x[1]-1)*n_y + id_y$

    $B = A + n_y$

    $C = B + n_y$

    $D = C + n_y$

    $p2N = [A,\; B,\; C,\; D]$

    \item[(4)] Chop off and convolute shape functions and derivatives in multi-dimensions through the tensor product operation:

    $N_{Ip} = \text{kron}({N}_{i,q}(x)[id_x], {N}_{i,q}(y)[id_y])$

    $\partial_x N_{Ip} = \text{kron}(\frac{\text{d}{N}_{i,q}(x)}{\text{d}x}[id_x], {N}_{i,q}(y)[id_y])$

    $\partial_y N_{Ip} = \text{kron}({N}_{i,q}(x)[id_x], \frac{\text{d}{N}_{i,q}(y)}{\text{d}y}[id_y])$

    $\nabla N_{Ip} = [\partial_x N_{Ip};\; \partial_y N_{Ip}]$
    
 \end{enumerate}

\section*{Nomenclature}
\label{nom}

\begin{itemize}
\item $(\bullet)_p$: a notation associated with material point
\item $(\bullet)_I$: a notation associated with grid node
\item $(\bullet)^e$: the elastic part of $(\bullet)$
\item $(\bullet)^p$: the plastic part of $(\bullet)$
\item $(\bullet)^v$: the volumetric part of $(\bullet)$
\item $(\bullet)^d$: the deviatoric part of $(\bullet)$
\item $(\bullet)^k$: $(\bullet)$ at $k^{th}$ number of Newton-Raphson iteration
\item $\delta$: the increment of $(\bullet)$ in the Newton-Raphson iteration
\item $\Delta$: the increment of $(\bullet)$
\item $\nabla$: the gradient of $(\bullet)${, which is the gradient in the undeformed configuration.}
\item $\nabla_x$: the spatial gradient of $(\bullet)$ {, which is the gradient in the deformed configuration.}
\item $\Bar{(\bullet)}$: the averaged $(\bullet)$
\item $\pmb{a}$: the acceleration
\item $\pmb{b}$: the left Cauchy-Green strain
\item $B$: the width of the footing
\item $c$: the cohesion
\item $c_f$: the factored cohesion
\item $\pmb{D}$: the consistent tangent matrix
\item $\text{d}t$: the time increment
\item $dof$: the spatial dimension of the 
\item $E$: Young’s modulus
\item $f$: the error of Newton-Raphson iteration
\item $\pmb{F}$: the deformation gradient
\item $\pmb{Fext}$: the external force
\item $\pmb{Fint}$: the internal force
\item $\pmb{Foob}$: the out of balance force
\item $\pmb{G}$: the strain-displacement matrix
\item $H$: the height of the slope
\item $\pmb{I}$: the identity matrix
\item $J$: the Jacobian of the deformation gradient
\item $k$: the number of Newton-Raphson iteration
\item $\pmb{K}$: the global stiffness matrix
\item $m$: the mass
\item $\pmb{M}$: the global mass matrix
\item $n$: the total number of B-spline shape functions
\item $n_n$: the number of grid nodes that influence the material point
\item $n_p$: the number of material points associated with node $I$
\item $N_{Ip}$: the shape function in the material point method
\item $q$: the polynomial order
\item $t$: the time
\item $tol$: the tolerance of Newton-Raphson iteration
\item $\pmb{u}$: the displacement
\item $\pmb{v}$: the velocity
\item $V$: the volume
\item $\pmb{x}$: the spatial coordinates
\item $x_i$: the $i^{th}$ component of knot vector
{\item $\pmb{X}$: the material coordinates}
\item $z$: the penetration depth of the footing
\item $\gamma$: the unit weight
\item $\delta_{ij}$: the Kronecker delta
\item $\delta_{max}$: the maximum displacement in the slope
\item $\pmb{\epsilon}$: the logarithmic strain
\item $\nu$: Poisson’s ratio
\item $\Xi$: the knot vector
\item $\pmb{\sigma}$: the Cauchy stress
\item $\pmb{\tau}$: the Kirchhoff stress
\item $\phi$: the friction angle
\item $\phi_f$: the factored friction angle
\item $\psi$: the dilation angle

\end{itemize}

\bibliography{lib}%

\end{document}